\documentclass[12pt]{article}
\pdfoutput=1
\usepackage{graphicx}
\usepackage{hyperref}
\usepackage{amsmath}
\usepackage{amssymb}

\usepackage[paperwidth=595pt,paperheight=841pt,top=70pt,right=70pt,bottom=70pt,left=70pt]{geometry}
\title{Generalization of Non-periodic Rhomb Substitution Tilings.}
\author{T.Hibma\\Zernike Institute for Advanced Materials\\University of Groningen\\t.hibma@rug.nl}

\begin{document}

\maketitle

\begin{abstract} General substitution rules for non-periodic rhomb tilings are derived. From the requirement that all substitution tiles consist of a discrete number of prototiles, it follows that a substitution tile with angle $s\pi{}/n$ must be built out of pairs of prototiles with angles $\left(s\pm{}t\right)\pi{}/n$, except if $t=0$. In addition, we require that a discrete number of prototile edges must fit between the beginning and endpoint of the substitution tile edge.  By comparing the total area of the discrete number of prototiles constituting the substitution tile and the total area derived from the assumed edge shape, a set of substitution rules is derived. The generalization involves the introduction of prototiles with a negative area or subtraction tiles.  
\end{abstract}

\section{Introduction}

The interest in non-periodic tiling is shared by seemingly different groups of people, such as crystallographers, mathematicians, artists, or just individuals interested in the beauty of geometric patterns. The discovery of the first quasicrystals \cite{Shechtman84} in the eighties of the last century triggered many solid state physicists and chemists to explore solids having structures intermediate between crystallographically ordered and completely amorphous materials. Mathematicians designed new concepts to distinguish and classify structures showing aperiodic order. An extensive description of the achievements up to now is given in the excellent book by Michael Baake and Uwe Grimm \cite{Baake13}. 

A powerful way to generate an aperiodic tiling is by the repetitive substitution of a set of basic tiles, the socalled prototiles.  The most commonly used \textit{tiling substitution rules} contain a prescription of the way in which the linearly expanded prototiles are replaced by a number of original prototiles. This type of substitution rule is called geometric in\cite{Priebe08}, to distinguish it from a combinatorial one. The most famous substitution tilings based on pairs of rhomb tiles with opening angles $s\pi{}/n, s=1,2$ are the Ammann-beenker ($n=4$) \cite{Beenker82} and the Penrose ($n=5$) \cite{Penrose74} \cite{Gardner77} tiling. For the Socolar ($n=6$) \cite{Socolar89} and Goodman-Strauss($n=7$) \cite{HarrissFrett}  tiling a set of three rhomb tiles are needed. The Goodman-Strauss tiling was later generalised to arbitrary $n>3$ by Harriss \cite{Harriss05}. His general substitution tile is composed of 16 prototiles, having edges with either a dent or a dimple in the middle with an edge angle of  $\pi/n$. Special substitution rules had to be designed for the smaller substitution tiles. In the following section we will show that an even more general approach is possible.

\section{A general model to construct rhomb substitution tiles}
\label{sec:GeneralModel}

The model is based on the observation that the combined area of the pair of rhomb prototiles $T_{s+t}^n$ and $T_{s-t}^n$ is proportional to the area of the rhomb prototile $T_s^n$, the proportionality factor $2\cos(t\pi/n)$ being independent of $s$. Consequently, substitution tiles $S_s^n$ may be constructed from a combination of prototiles $T_s^n$ and pairs of tiles $T_{s\pm t}^n$. If $n_0$ and $n_t$ are their numbers respectively, the area of the substitution tile is 

\begin{equation}
  S=n_0+\Sigma_t 2n_t\cos(t\pi/n)
  \label{eq:AreaLeft}
\end{equation}

\begin{figure}[!ht]
	  \centering
		\includegraphics[width=0.40\textwidth]{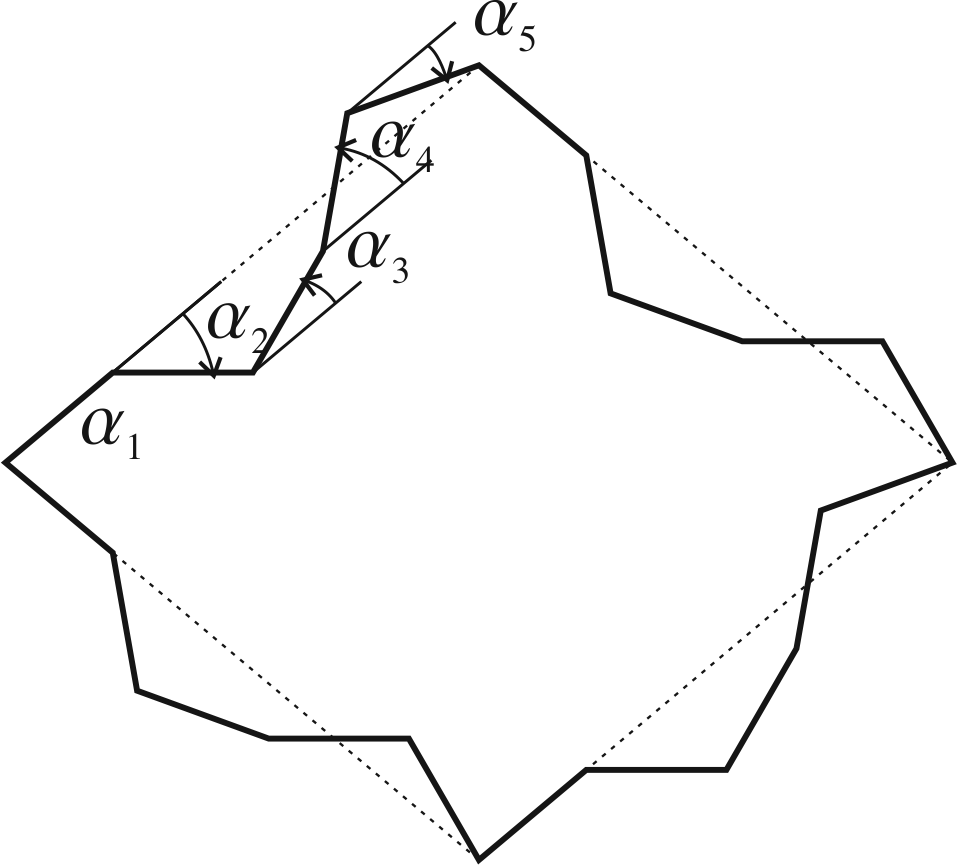}
	  \caption{\small { Substitution tile edge structure. The edge angles $\alpha_i$ occur in $\pm$ pairs or are zero. In this example we chose the edge angles $\alpha_i$ to be $0$, $-2\pi/9$ ,  $\pi/9$,  $2\pi/9$, $-\pi/9$ respectively. All edges have congruent shapes. The lower and upper left edges are related by a rotation over the opening angle $s\pi/9$ with respect to the left corner. Similarly, the lower and upper right edges are related by a rotation over the opening angle with respect to the right corner. Opposite edges are related by a translation.}}
	  \label{fig:SubstitutionTileEdgeShape}
\end{figure}

A second requirement for a tiling of the entire plane is to realize proper edge substitutions. We will assume that all four substitution tile edges have the same shape. Neighboring edges at the opening angle are related by a rotation over that angle and opposite edges are related by a translation (figure \ref{fig:SubstitutionTileEdgeShape}). This edge arrangement also ensures that the substitution tile area is equal to the inflated rhomb area, and, therefore, $S$ is equal to the \textit{areal scaling factor}.  The angles between the outer prototile edges and the substitution tile rhomb edge will be called the \textit{edge angles} $\alpha_i$. For now, we will assume that overhangs are not allowed and $|\alpha_i|\leq \pi{}/2$. The length of the substitution tile rhomb edge

 \begin{equation}
  L=\Sigma_i\cos\alpha_i
  \label{eq:InflationFactor}
\end{equation}  

is the \textit{inflation factor} of the rhomb tiles. Because the areal scaling factor $S$ is the square of $L$, equations \ref{eq:AreaLeft} and \ref{eq:InflationFactor} can be combined into

\begin{equation}
  n_0+\Sigma_t 2 n_t\cos(t\pi /n)=\frac{1}{2}\Sigma_i\Sigma_j\{\cos(\alpha_i+\alpha_j)+\cos(\alpha_i-\alpha_j )\} 
  \label{eq:AreaComp}
\end{equation} 
                     
This equality can only be satisfied, if the arguments $\alpha_i+\alpha_j$ and $\alpha_i-\alpha_j$ are both equal to an integer times $\pi/n$ for all $i$ and $j$. There are two solutions: either all angles $\alpha_i$ are equal to an integer times $\pi/n$, or all of them are equal to a half-integer times $\pi/n$. 

Because the beginning and end of the substitution edge have to be at the endpoints of the substitution tile rhomb edge, the following relationship between the edge angles $\alpha_i$ should be met: 
\begin{equation}  
  \Sigma_i sin\alpha_i =0 
  \label{eq:SumEdgeSines}
\end{equation} 
                                                        
A general solution is that the edge angles $\alpha_i$ occur in $\pm$-pairs or are zero 
\footnote{There are also special solutions. For instance, if one requires that the sum of three terms are zero, one finds that $\alpha_1+\alpha_2=\pi{}/3$ and $\alpha_1+\alpha_3=-\pi{}/3$. This solution is valid if $n$ is a multiple of 3. An example satisfying this condition is the Lord tiling, having edge angles $\alpha_1=\alpha_2=\pi{}/6$ and $\alpha_3=-\pi{}/2$ \cite{HarrissFrett}. $n=3$ in this case, and the edge sequence is $(\frac{1}{2},-\frac{3}{2},\frac{1}{2})$. In this paper, however, we will only consider the more general $\pm$ pairing condition.}
, and equation \ref{eq:AreaComp} becomes 

\begin{equation}  
n_0+2\Sigma_t n_t \cos(t\pi/n) =\{m_0+2\Sigma_r m_r  \cos(r\pi/n)\}^2  
\label{eq:AreaInt}
\end{equation}    
                                           
or
\begin{equation}  
 n_0+2\Sigma_t n_t\cos(t\pi/n)=\{2\Sigma_r m_{r+\frac{1}{2}}\cos((r+\frac{1}{2})\pi/n)\}^2 
\label{eq:AreaHalfInt}
\end{equation}

Equations \ref{eq:AreaInt} or \ref{eq:AreaHalfInt} determine the type and number of prototiles $T_{s\pm t}^n$ from which the substitution tile can be constructed, once the shape of the substitution tile edge has been chosen.  In view of the above considerations, this edge shape may be characterized by a sequence of integers or halfintegers, the \textit{edge sequence} $(k_i)_1^N$, defined by $\alpha_i=k_i\pi{}/n$, $-n/2\leq k_i\leq n/2$ \cite{Maloney14}.

If the finite edge angles are present as pairs in accordance with equation \ref{eq:SumEdgeSines}, always a valid solution for the substitution tile is obtained, because both sides may be written as a sum of cosine terms having even valued coefficients. The pairing of the edge angles, therefore, guarantees that the substitution tiles are composed of an integer number of prototiles.

\begin{figure}[!ht]
	\centering
		\includegraphics[width=0.40\textwidth]{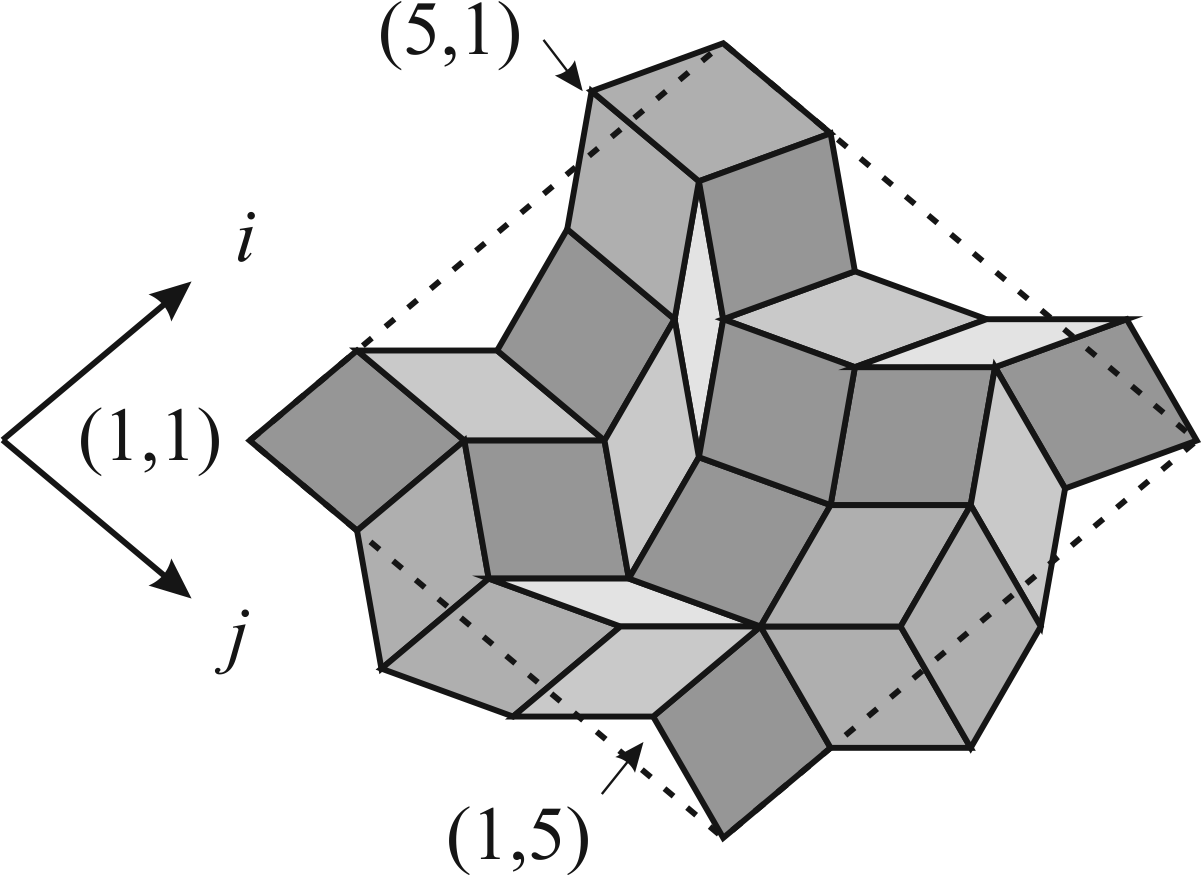}
	\caption{\small {General substitution scheme for a substitution tile  $S_s^n$ with a given edge shape. The prototile at position $i$,$j$  is  $T_{s+k_i-k_j}$, where $k_i\pi/n$  and $k_j \pi/n$ are the edge angles at the upper and lower left edges respectively.}}
	\label{fig:GeneralScheme}
\end{figure}

Equations \ref{eq:AreaInt} or \ref{eq:AreaHalfInt} constitute a connection between the prototile edge angle pairs $m_k$ and the numbers of prototiles in a substitution tile $n_t$, not their arrangement.  The relations do not guarantee that a consistent set of substitution tiles can be found. However, in the following we will show that a general set of substitution rhomb tiles can be constructed for arbitrary $n$ and for an arbitrary substitution tile edge shape. 

We start with a construction of the circumpherence of the tile $S_s^n$ as described earlier and illustrated in figure \ref{fig:SubstitutionTileEdgeShape}. Next, copies of the edges are translated to the breaks of neighboring edges. If the breaks of the upper and lower left edge are indexed as $i\in(1,2..n)$ and $j\in(1,2..n)$ respectively, starting at the left corner as indicated in figure \ref{fig:GeneralScheme}), one obtaines a grid of vertices $(i,j)$, at which four prototiles meet. The one bounded by the vertices $(i,j)$, $(i+1,j)$, $(i,j+1)$ and $(i+1,j+1)$ is a prototile of the type $T_{s+k_i-k_j}^n$. The vertices at diagonal positions are occupied by tiles $T_s^n$ , whereas one can find pairs of tiles $T_{s\pm(k_i-k_j )}^n$ at off-diagonal positions $(i,j)$ and $(j,i)$. This general substitution rule may be represented by the following matrix 

\begin{equation}
S_s^n=
\begin{bmatrix}
T_s^n &T_{s+k_1-k_0}^n & T_{s+k_2-k_0}^n & \ldots &T_{s+k_{n-1}-k_0}^n  \\
T_{s+k_0-k_1}^n &T_s^n & T_{s+k_2-k_1}^n & \ldots & T_{s+k_{n-1}-k_1}^n  \\
T_{s+k_0-k_2}^n &T_{s+k_1-k_2}^n & T_s^n & \ldots & T_{s+k_{n-1}-k_2}^n  \\
\vdots & \vdots & \vdots & \ddots & \vdots \\
T_{s+k_0-k_{n-1}}^n &T_{s+k_1-k_{n-1}}^n &T_{s+k_2-k_{n-1}}^n & \ldots & T_s^n \\
 \end{bmatrix}
\label{eq:Substmatrix} 
\end{equation}

For later use one should note, that the prototiles parallel to the substitution edges, i.e. the rows or columns of the matrix, form worms, and the edges of the worms have shapes identical to the edge shape of the substitution tile. 

The prototiles are allowed to have indices $s\pm t<0$ or $>n$. These prototiles will have negative areas, meaning that they have to be subtracted from the tiling. We consider a tiling of the plane to be a legitimate one, if in the end there are no holes or overlaps. So, negative or subtraction tiles are allowed, if they remove all overlaps between tiles and do not leave holes in the tiling.  In one of the next sections, we will reason, that this is presumably the case for substitution edges without loops. Also the zero area prototiles for which $s\pm t=0$ or $n$ play a important role in our scheme and cannot simply be neglected.

\section{Substitution tiles with $m_1=1$}

We will now discuss one of the simplest examples of this general scheme, i.e. the one with $m_1=1$ and $m_{>1}=0$.  From \ref{eq:AreaInt} it follows that
\begin{equation}  
 n_0=m_0^2+2m_1^2 ,\  n_1=2m_0 m_1 ,\  n_2=m_1^2
\end{equation}
For $m_0=0$ we get  $n_0=2$ and $n_1=1$  and the $2\times2$ general substitution tile is
\begin{equation}
S_s^n=
\begin{bmatrix}
T_s^n &T_{s+2}^n \\
T_{s-2}^n &T_s^n \\
 \end{bmatrix}
\label{eq:Sm1=1matrix} 
\end{equation}
The edge sequence is (1, -1) and the inflation factor $L=2 \cos(\pi/n)$.
                                                         
\begin{figure}[!ht]
	\centering
		\includegraphics[width=0.9\textwidth]{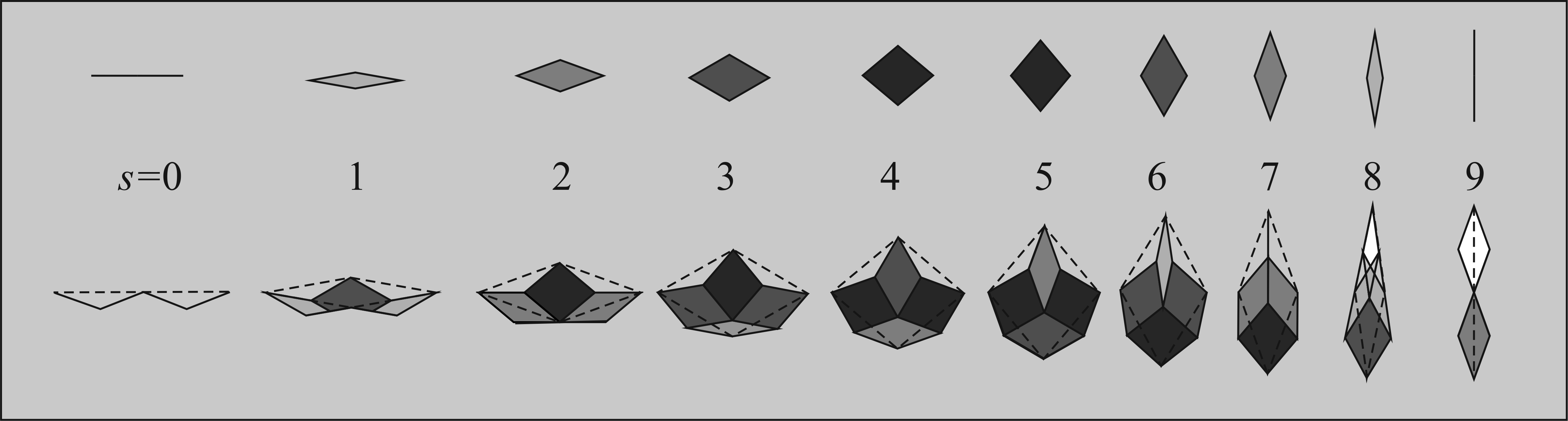}
	\caption{\small $m_1=1$ tiles for $n=9$. The $s=0,1,8$ and $9$ tiles have been constructed using negative area prototiles.}
	\label{fig:Full-dent-9-2x2-tiles}
\end{figure}

In figure \ref{fig:Full-dent-9-2x2-tiles} all first generation substitution tiles, constructed using the general scheme outlined in the previous section, are shown for $n=9$. Negative areas are colored white in a gray background. The $s=0$, $1$, $8$ and $9$ substitution tiles are constructed in exactly the same way as the general tile, i.e. the $s$ tiles are rotated over $\pm\pi/n$ and connected with the unrotated $s\pm2$ tiles. For $s=0$ two zero area prototiles $T_0^n$, which are lines of length 2, a $T_2^n$ prototile and a $T_{-2}^n$ subtraction prototile have to be used. The normal and subtraction tile annihilate each other, and the substitution tile is reduced to a zigzag line. This is most easily seen by considering tiles $T_{\delta}^n$, and $T_{\delta \pm 2}^n$ letting $\delta$ approach $0$. For $s=n$, on the other hand, the $T_{n+2}^n$ and the $T_{n-2}^n$ tiles do not overlap and the substitution tile consists of a positive and negative part. However, the total area is zero, as it should. For $s=1$ the negative area $T_{-1}^n$ tile cuts pieces off the other three prototiles. For $s=n-1$ the two prototiles $T_{n-1}^n$ partly overlap. This overlap is compensated  by the negative area $T_{n+1}^n$ tile, leaving a fractional negative region shown in white. An essential result is that all these tiles have the same dents and dimples as the general tile and can be used to construct larger substitution tiles.

\begin{figure}[!ht]
	\centering
		\includegraphics[width=0.9\textwidth]{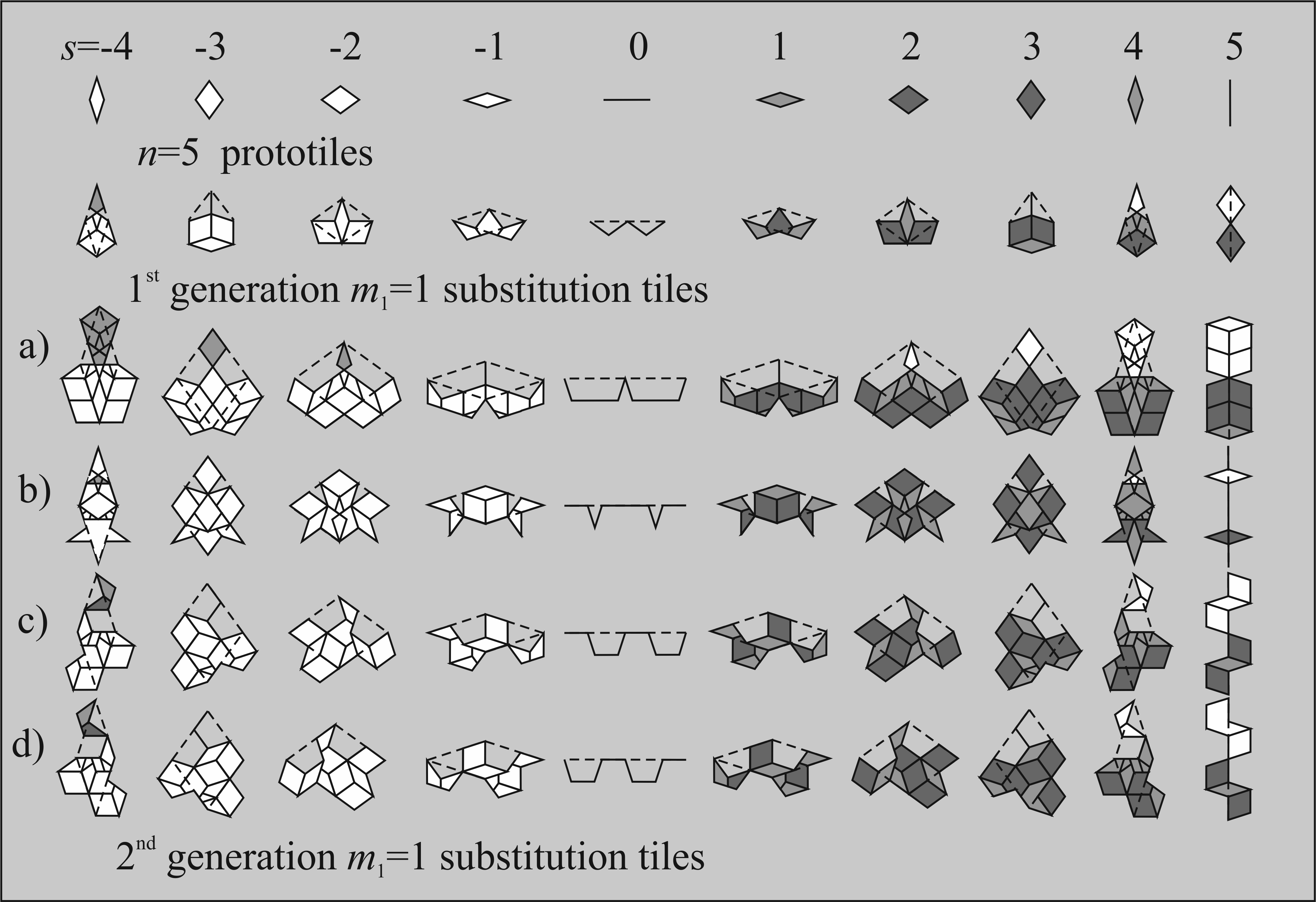}
	\caption{\small {$m_1=1$  tiles for $n=5$.  White parts have a negative value and are to be subtracted from the tiling. Four different sets of substitution tiles can be constructed for the second generation. }}
	\label{fig:Fulldent-5-2x2-tiles}
\end{figure}

\begin{figure}[!ht]
	\centering
		\includegraphics[width=0.9\textwidth]{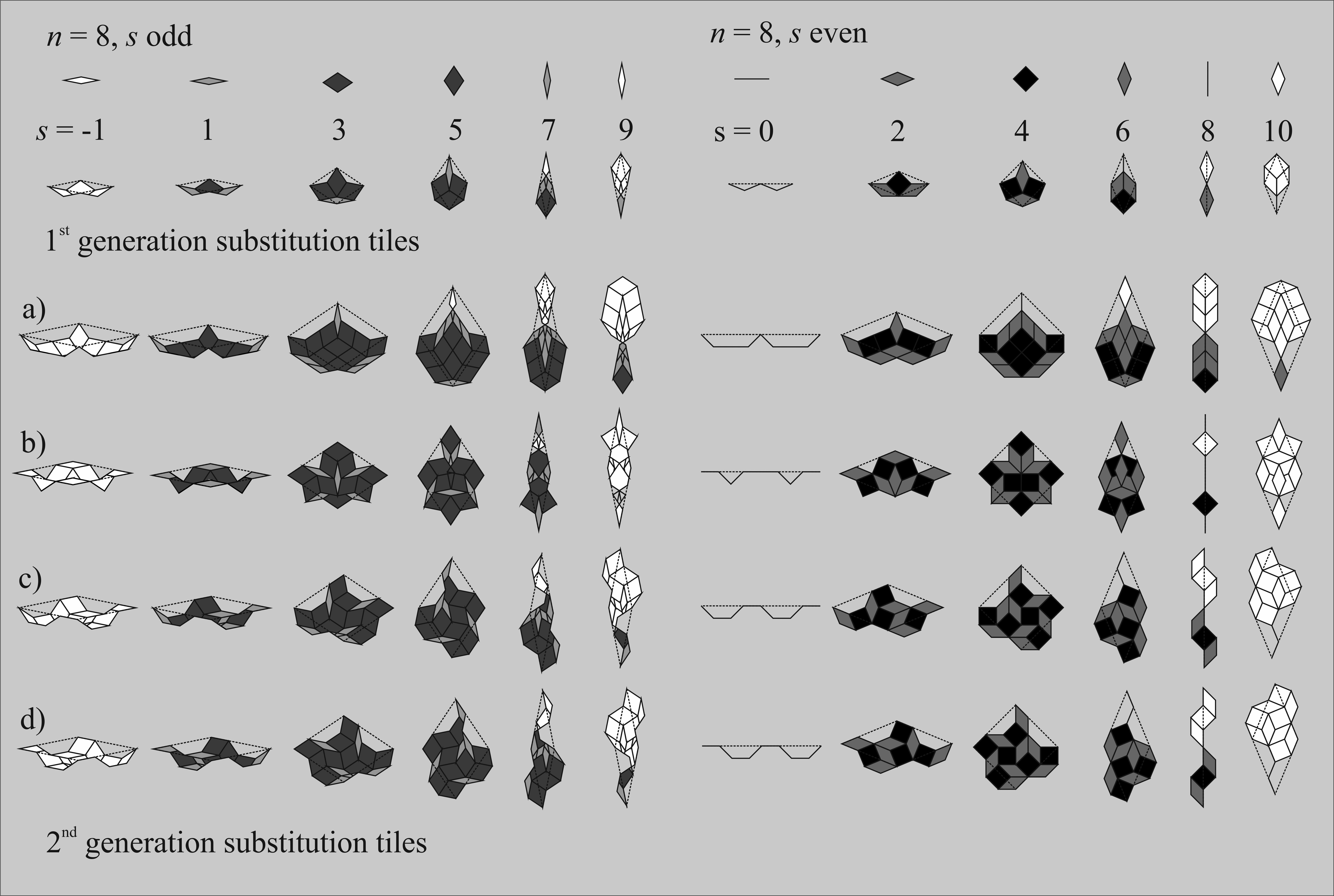}
	\caption{\small {$m_1=1$  tiles for $n=8$. Odd and even $s$ tiles form separate tile subsets, which do not mix.The even $s$ subset is identical to the $m_\frac{1}{2}$, $n=4$  tile set.}}
	\label{fig:Fulldent-8-2x2}
\end{figure}

\begin{figure}[!ht]
	\centering
		  \includegraphics[width=0.8\textwidth]{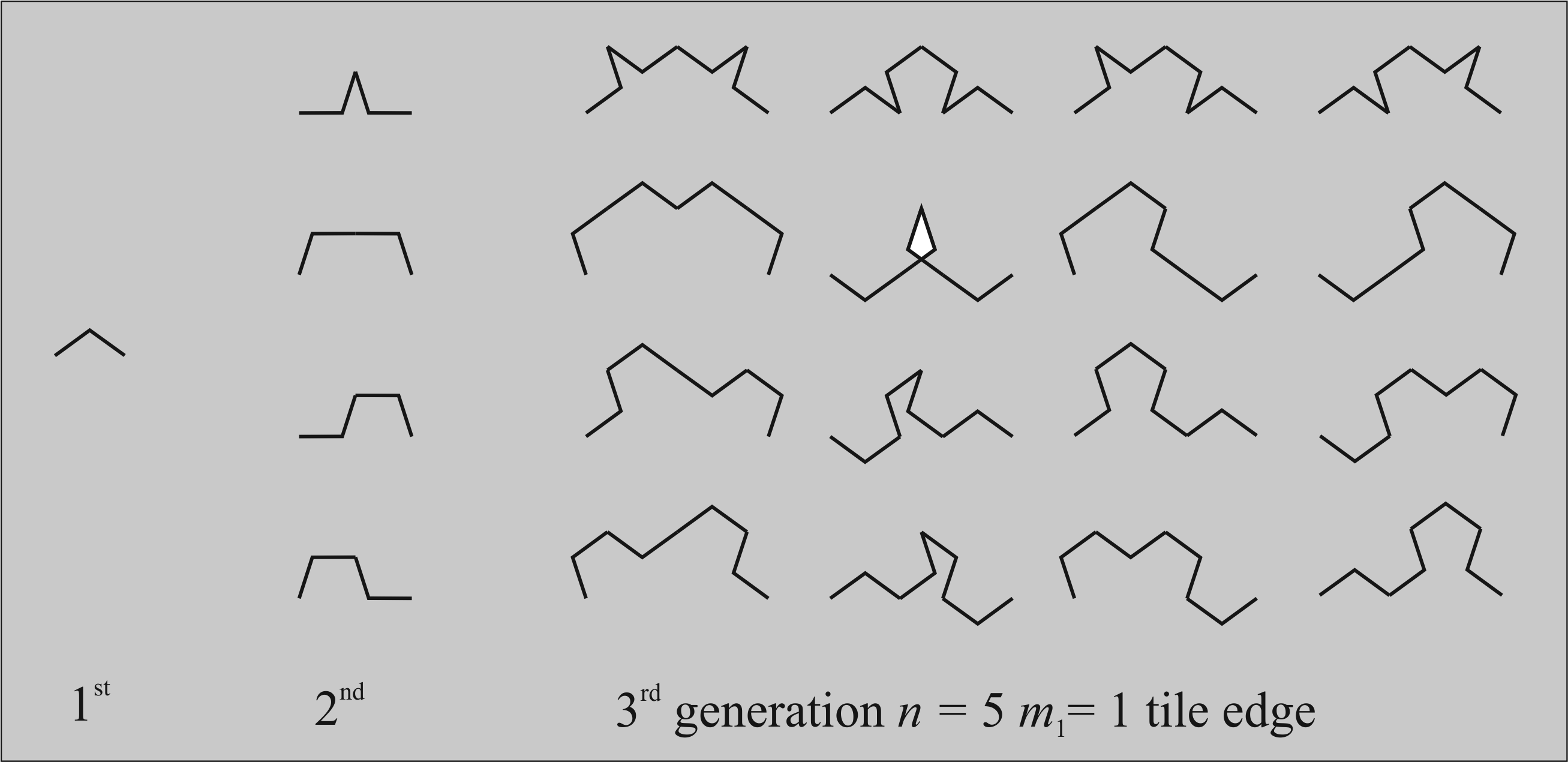}
	    \caption{\small {Edge shapes for three generations of $n=5$ $2\times 2$ full dent substitution tiles. }}
      \vspace{10pt}
	    \label{fig:Harris-5-2x2-edges}
		   \includegraphics[width=0.8\textwidth]{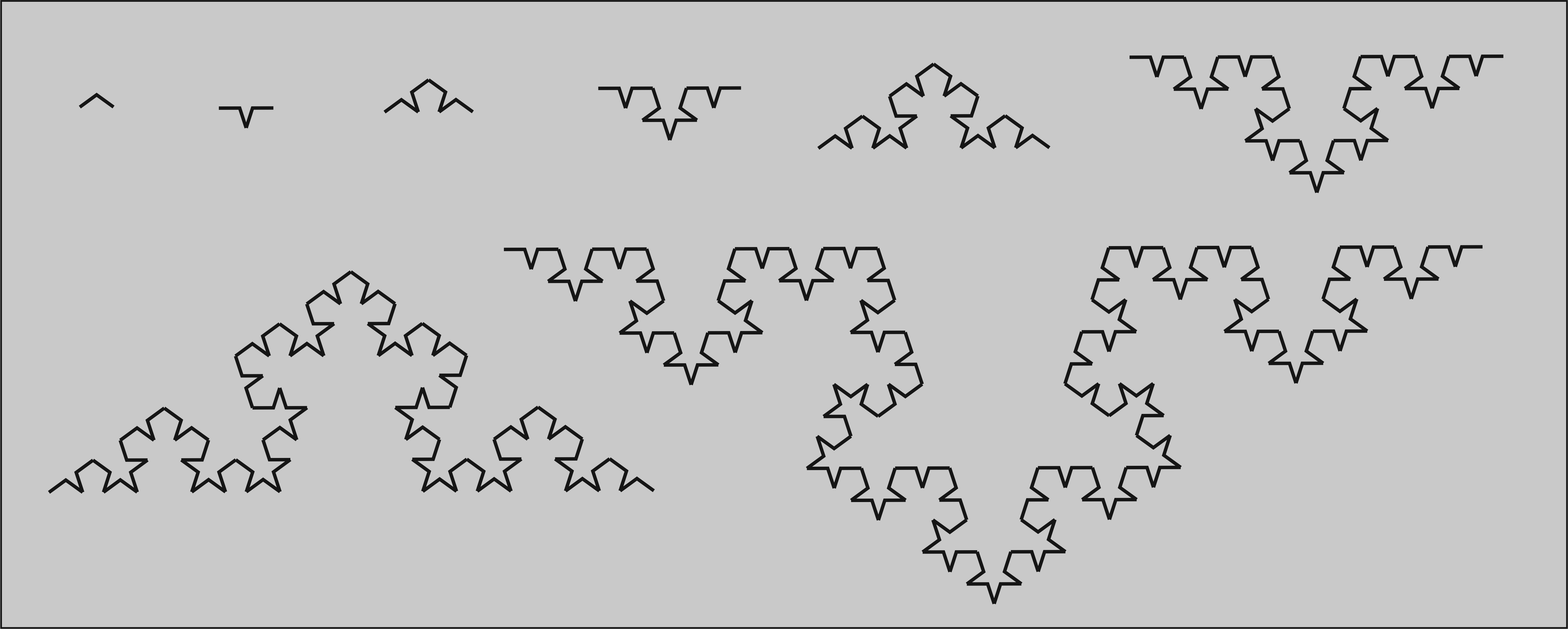}
	    \caption{\small {Example of substitution tile edge development leading to Koch tile edge.}}
      \vspace{10pt}
	    \label{fig:KochEdges}
		   \includegraphics[width=0.8\textwidth]{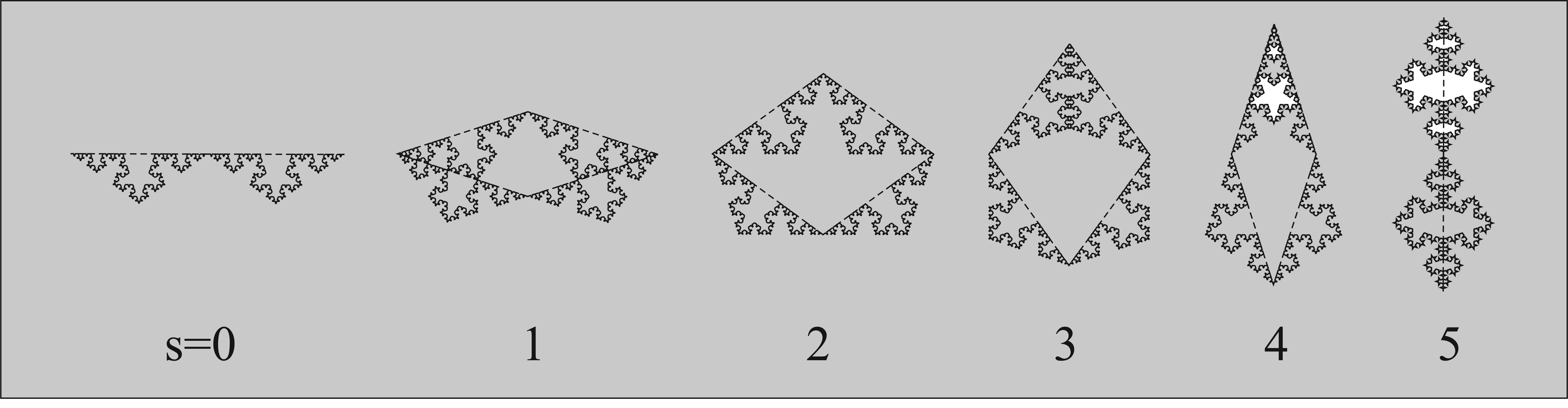}
	    \caption{\small {Aperiodic $n=5$ von Koch Rhomb Tiles. }}
      \vspace{10pt}
	    \label{fig:KochTiles}
\end{figure}

The prototiles  $T_s^n$ and $T_{n-s}^n$ have identical rhomb shapes and, therefore, one obtains an identical object if one replaces the other. The same holds if a prototile is rotated over $\pi$.  However, this is different for higher generation substitution tiles. Although the substitution tiles $S_s^n$ and $S_{n-s}^n$ have the same basic rhomb shapes, their edge configurations are different. A replacement or a twofold rotation of one or more of the four prototiles, leads to 16 possible tiles with different edge shapes. The reasoning is as follows: if we fill up a next generation $S_s^n$ tile, there are four ways to put a $T_s^n$ or $T_{n-s}^n$ tile at the ($1$, $1$) position. Subsequently, two out of four $T_{s\pm 2}^n$ or $T_{n-s\mp 2}^n$ tiles fit at the ($1$, $2$) and ($2$, $1$) position. Finally, only one out of four $T_s^n$ or $T_{n-s}^n$ type tiles matches the edges already present at the ($2$, $2$) position. The total number of tile shapes is, therefore, sixteen. However, only four of them are consistent with the requirement that the four edge shapes are the same and related to each other by a translation or rotation as defined in figure \ref{fig:SubstitutionTileEdgeShape}. These are: a) the basic substitution tile, b) a substitution tile in which all constituent tiles in a) have been rotated over $\pi$, c) a substitution tile in which the $s\pm 2$ tiles are replaced by $n-s\mp 2$ tiles and d) a substitution tile in which all constituent tiles in c) have been rotated over $\pi$. Another, simpler way to arrive at the four substitution rules is to consider the edge sequence of a tile. There are two ways to connect an edge between two corners of a rhomb. Consequently, the number of different edge shapes for a substitution tile is $2^N$, if $N$ is the number of prototile edges along the substitution tile edge, which is two in the case considered here. The four second generation substitution tile varieties are shown for $n=5$ in figure \ref{fig:Fulldent-5-2x2-tiles} and $n=8$ in figure \ref{fig:Fulldent-8-2x2}. The first two are mirror symmetric, whereas the last two are asymmetric, but mirror images of each other. All four sets contain tiles with negative area parts. These negative areas account for overlapping dimples of neighbouring edges. They will naturally disappear if the next larger substitution tile is constructed, i.e. overlapping dimples will be filled up with dents of neighboring tiles. In some cases, notably the $s=n$ and $s=n-1$ tiles, new negative areas will appear at one end of the tile. Also in most of the basic substitution tiles of set a), the upper part is a negative area part. Second generation substitution tiles with $0<s<n-1$ do not have negative area parts for sets b), c) and d). In set b) this is due to the fact, that the $s=n$ and $s=n-1$ constituent tiles are used in a rotated position, so that no negative area parts stick out. In sets c) and d) the $s=n$ and $s=n-1$ prototiles have been replaced by the $s=0$ and $s=1$ prototiles respectively, not having negative area parts by themselves.

Because the $s$-values of the prototiles in substitution tiles belonging to set a) and b) differ by 2, the even and odd $s$ substitution tiles constitute separate subsets, both fully tiling the plane. For even $n$ this is also true for sets c) and d), whereas for odd $n$ the even and odd tiles intermix.    

Some of the prototiles are cut into smaller and rather complicated shapes by the subtraction procedure. The latter is not unique, because the final surface is obtained by cutting the shape of the negative area out of one of two overlapping dents.  In figure \ref{fig:Fulldent-5-2x2-tiles} this region has been filled in quite arbitrarily by showing a linepattern consistent with the symmetry of the tile. Although we cannot prove it rigorously, it appears that it will be possible to replace the mutilated tiles by half-prototiles (triangles) for $n$=5. 

Each new generation of $m_1=1$ substitution tiles will be fourfold in the same way as was described for the second generation. For each new generation, an edge segment is replaced by a two-segment dent or dimpel. Already after a few generations a large number of different edge shapes are generated. In figure \ref{fig:Harris-5-2x2-edges} the edge shape is shown for three generations for $n=5$. Sometimes, overhangs or loops show up. The latter are due to overlapping neighboring dimples causing negative area parts at the circumpherence of a substitution tile.  In figure \ref {fig:Harris-5-2x2-edges} this is seen to occur once. The way a next generation edge shape is obtained is reminiscent of the construction of a von Koch curve \cite{vonKoch1904}. An important difference is that in the original von Koch construction each segment is replaced by four instead of two new segments. The actual von Koch curve is obtained if substitution rule b) of figure \ref {fig:Harris-5-2x2-edges} is applied for each new generation. Eight generations of  these von Koch edges are shown in figure \ref{fig:KochEdges}, and the corresponding von Koch tiles in figure \ref{fig:KochTiles}. A complete tiling of the plane may be generated using either the odd or even $s$ tiles.

\section{Substitution tiles with $m_\frac{1}{2}=1$}
\label{Sec:FracMTiles}

In section \ref{sec:GeneralModel} we found that there are two types of edge substitution rules, one with integer and one with half-integer edge angle fractions of $\pi$, eqs. \ref{eq:AreaInt} and \ref{eq:AreaHalfInt} respectively. In this section we will discuss the second case. From equation \ref{eq:AreaHalfInt}

\begin{equation} 
n_0=2m_\frac{1}{2}^2 ;\   n_1=m_\frac{1}{2}^2 ;\   n_{>\frac{1}{2}}=0   
\end{equation}

\begin{figure}[!ht]
	\centering
		\includegraphics[width=0.9\textwidth]{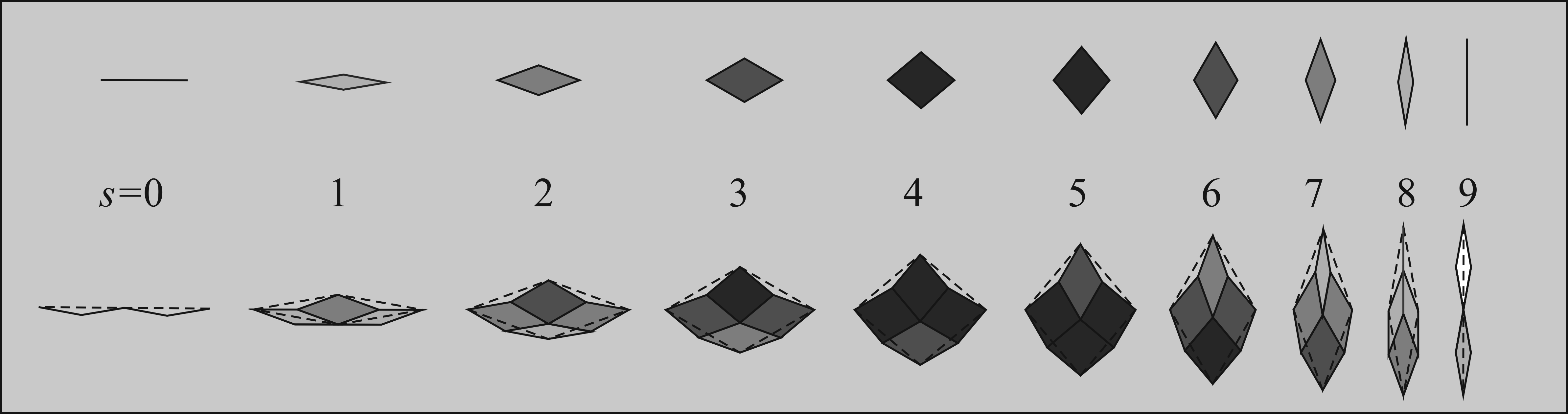}
	\caption{\small {$m_\frac{1}{2}=1$ substitution tiles for $n=9$.}}
	\label{fig:Halfdent-9-2x2-tiles}
\end{figure}

\begin{figure}[!ht]
	\centering
		\includegraphics[width=0.9\textwidth]{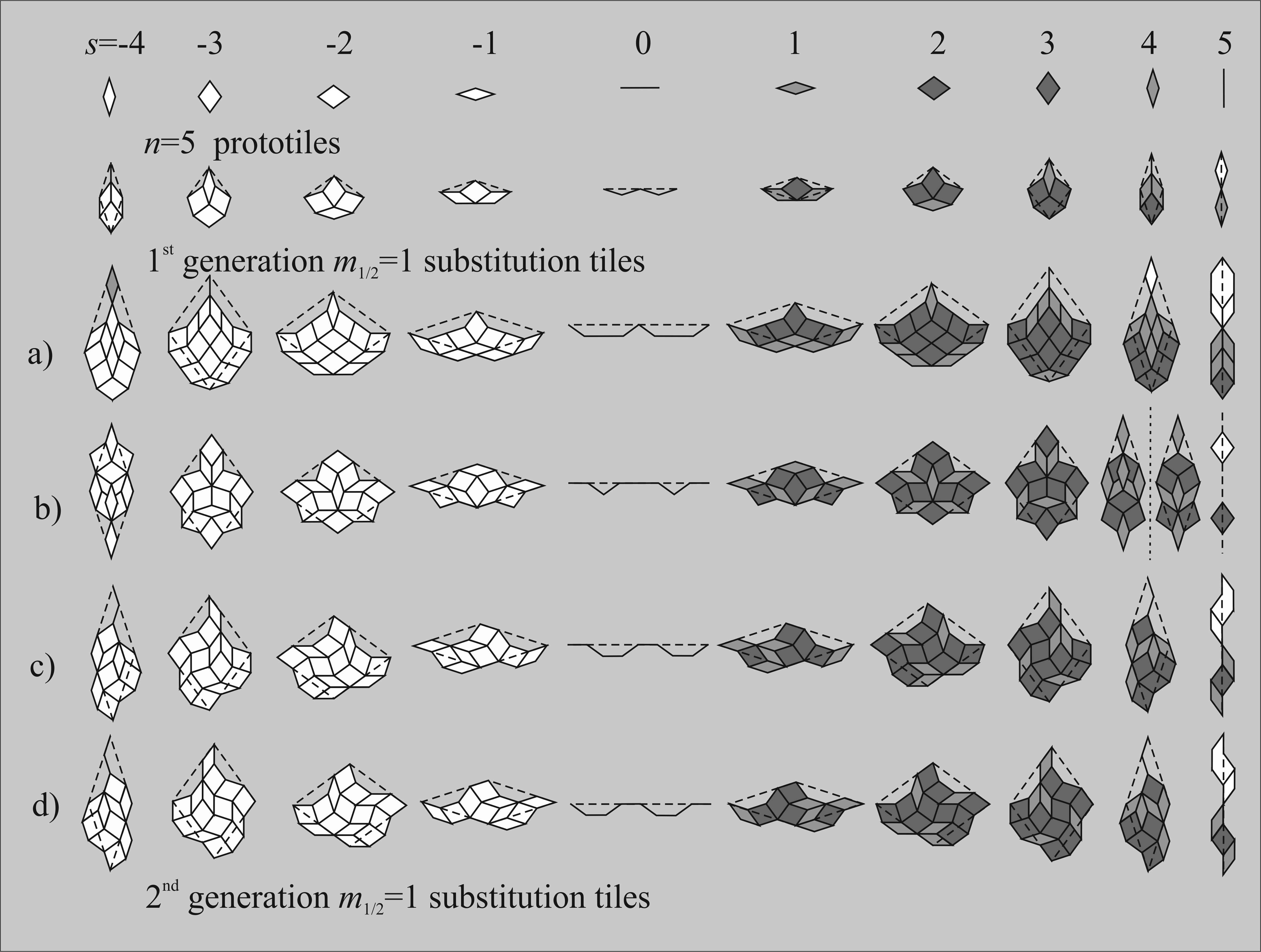}
	\caption{\small {$m_\frac{1}{2}=1$ $n=5$ (Penrose Rhomb) tiles. White denotes a negative area. Like in the $m_1=1$ case, four different 2nd generation tile sets may be constructed, two mirror symmetric sets a) and b) and two asymmetric sets c) and d). The edge sequence of the second symmetric set (b) is identical to the one of the Harriss $n=5$ set of tiles. The odd s subset of c) and d) is the so-called binary tiling \cite{HarrissFrett}}}
	\label{fig:Halfdent-5-2x2-tiles}
\end{figure}  
                                              
If $m_\frac{1}{2}=1$, it follows that $n_0=2$ and $n_1=1$. Thus, like in the case of $m_1=1$, four prototiles are involved, but the non-congruent tiles have $s$-values differing only by one instead of two, or

\begin{equation}
S_s^n=
\begin{bmatrix}
T_s^n &T_{s+1}^n \\
T_{s-1}^n &T_s^n \\
 \end{bmatrix}
\label{eq:Sm1/2=1matrix}
\end{equation}
The edge sequence is $(\frac{1}{2}, -\frac{1}{2})$ and the inflation factor $L=\cos \pi{}/2n$.
 
In figure \ref{fig:Halfdent-9-2x2-tiles} the 1st generation substitution tiles for $n=9$ are shown. In contrast to the $m_1=1$  substitution tiles (figure \ref{fig:Full-dent-9-2x2-tiles}), negative area tiles are not needed for the construction of positive area tiles. Also in this case, there are four possible edge shapes for the second generation of substitution tiles, depending on how the first generation substitution tiles are combined (see figure \ref {fig:Halfdent-5-2x2-tiles} for $m_\frac{1}{2}=1$, $n=5$ tilings). The $s=n-1$ substitution tile of set (a) contains a negative area, because the dimples at neighboring edges overlap. In higher generation substitution tiles this negative area will also appear in the tiles with smaller $s$-values. For set b) the negative half of the $s=n$ constituent tile does not show up, because it removes the overlap of two $s=n-1$ tiles. Also in the higher generation tiles negative area will not be visible. In fact, the shape of the edge of the tiles of set b) is identical to the shape of the Harriss tiles. The tiles of set c) and d) are, like in the $m_1=1$ case, mirror images of each other. For odd $n$, they can be split up into two sets, both giving a full aperiodic tiling of the plane, one with odd and one with even $s$-values. The one with even $s$ tiles is special because no $s=n$ constituent tiles are needed and, therefore, no negative tiles are involved. By neglecting the the $s=0$ tiles, the set of prototiles is reduced to a conventional one. For $n=5$ the tiling based on the even $s$ set, i.e. $s=2$ and $s=4$ is already known as the "binary" or Lan\c{c}on Billard (LB) tiling \cite{Lancon88}. The odd $s$ set contains the $s=n$ tile. Its role is to remove a patch of prototiles at one side of an $s=1$ tile and add it to the other side. A tiling results which is very similar to the LB tiling. The prototile set for even $n$ comprises all $s$ tiles from 0 to $n$. Although the negative part of the $s=n$ tile is involved in the substitution rule, the prototiles are not cut up and a tiling consisting of complete prototiles results.

The $m_\frac{1}{2}=1$ substitution tiles are identical to a subset of $m_1=1$ tiles with a doubled $n$-value and with even $s$-values. So, the tiles shown in figure \ref{fig:Halfdent-5-2x2-tiles} are also obtained as an even $s$ subset of the  $m_1=1$, $n=10$ case and the even $s$ subset of the  $m_1=1$, $n=8$ tiles shown in figure \ref{fig:Fulldent-8-2x2} is identical to the $m_\frac{1}{2}=1$, $n=4$  tile set.

\section {Substitution Matrices.}

In this section we want to reformulate the rhomb substitution model in terms of the edge and tile substitution matrices.
If overhangs are included, the prototile edges in a tiling will point into $2n$ directions. Each of these is replaced by a number of prototile edges in orientations determined by the edge sequence, i.e. $m_0$ in the same, $m_n$ in the opposite direction and $m_k$ with $k\in (1, 2,.., n-1)$ in directions differing by $+k\pi{}/n$ and $-k\pi{}/n$. The edge substitution matrix, therefore, is

\begin{equation}
\textbf{M}=
\begin{bmatrix}
m_0 & m_1 & m_2 & \ldots & m_{n-1} & m_n  & \ldots & m_3 & m_2 & m_1  \\
m_1 & m_0 & m_1 & \ldots & m_{n-2} & m_{n-1} & \ldots &m_4 & m_3 & m_2  \\
m_2 & m_1 & m_0 & \ldots & m_{n-3} & m_{n-2} & \ldots &m_5 & m_4 & m_3  \\
\vdots & \vdots & \vdots & \ddots & \vdots & \vdots & \ddots  & \vdots & \vdots & \vdots \\
 m_{n-1} &  m_{n-2} &  m_{n-3} & \ldots &m_0 & m_1 & \ldots  & m_{n-2} & m_{n-1} &  m_n \\
 m_n & m_{n-1} &  m_{n-2}& \ldots & m_1 & m_0 & \ldots  &  m_{n-3} &  m_{n-2} &  m_{n-1} \\
\vdots & \vdots & \vdots & \ddots & \vdots & \vdots & \ddots  & \vdots & \vdots & \vdots \\
m_3 & m_4 & m_5 & \ldots  & m_{n-2} & m_{n-3} & \ldots& m_0 & m_1 & m_2  \\
m_2 & m_3 & m_4 & \ldots & m_{n-1} & m_{n-2} & \ldots & m_1 & m_0 & m_1  \\
m_1 & m_2 & m_3 & \ldots & m_n & m_{n-1} & \ldots & m_2 & m_1 & m_0  \\
\end{bmatrix}
\label{eq:Mmatrix} 
\end{equation}

A tile with index $s$ is substituted by $n_0$ prototiles with index $s$, $n_n$ with index $s+n$ and $n_t$ with index $s+t$ and $s-t$, with $t \in (1,...,n-1)$ and $s \in(0,...,2n-1)$.  So the substitution matrix $\textbf{S}$ is
 
\begin{equation}
\textbf{S}=
\begin{bmatrix}
n_0 & n_1 & n_2 & \ldots & n_{n-1} & n_n  & \ldots & n_3 & n_2 & n_1  \\
n_1 & n_0 & n_1 & \ldots & n_{n-2} & n_{n-1} & \ldots &n_4 & n_3 & n_2  \\
n_2 & n_1 & n_0 & \ldots & n_{n-3} & n_{n-2} & \ldots &n_5 & n_4 & n_3  \\
\vdots & \vdots & \vdots & \ddots & \vdots & \vdots & \ddots  & \vdots & \vdots & \vdots \\
 n_{n-1} &  n_{n-2} &  n_{n-3} & \ldots &n_0 & n_1 & \ldots  & n_{n-2} & n_{n-1} &  n_n \\
 n_n & n_{n-1} &  n_{n-2}& \ldots & n_1 & n_0 & \ldots  &  n_{n-3} &  n_{n-2} &  n_{n-1} \\
\vdots & \vdots & \vdots & \ddots & \vdots & \vdots & \ddots  & \vdots & \vdots & \vdots \\
n_3 & n_4 & n_5 & \ldots  & n_{n-2} & n_{n-3} & \ldots& n_0 & n_1 & n_2  \\
n_2 & n_3 & n_4 & \ldots & n_{n-1} & n_{n-2} & \ldots & n_1 & n_0 & n_1  \\
n_1 & n_2 & n_3 & \ldots & n_n & n_{n-1} & \ldots & n_2 & n_1 & n_0  \\
\end{bmatrix}
\label{eq:Nmatrix}
\end{equation}

The relation between the tile and edge substitution matrices is

 \begin{equation}
 \textbf{S}=\textbf{M}^2
\label{eq:NvsM} 
\end{equation}

This matrix equation may be used to calculate the numbers of prototiles in a substitution tile for a given edge shape instead of equations \ref{eq:AreaInt} . 
The $n_j$ are given by the product of the first row and the $j$-th column

\begin{equation}
n_j=\textbf{S}_{1 j+1}=\Sigma_{k=1}^{2n} \textbf{M}_{1 k}\textbf{M}_{k j+1}      ;j\in(0, 1,...., n) 
\end{equation}

Both $\textbf{M}$ and $\textbf{S}$ are \textit{circulant} $2n\times 2n$ matrices \cite{Kra12}. 
Consequently, a shorthand notation of equations \ref{eq:Mmatrix} and \ref{eq:Nmatrix} is

\begin{equation}
\textbf{M}=\textit{circ}(m_0, m_1, m_2, ..., m_{n-1}, m_n,...., m_3, m_2, m_1)
\end{equation}

\begin{equation}
\textbf{S}=\textit{circ}(n_0, n_1, n_2, ..., n_{n-1}, n_n,...., n_3, n_2, n_1)
\end{equation} 
All \textit{circulant} matrices are known to have the same set of normalized eigenvectors

\begin{equation}
 v_l={1\over\sqrt(2n)}(1, \epsilon^l, \epsilon^{2l}, \epsilon^{3l},.......,\epsilon^{(2n-1)l})^T
\end{equation} 
with $\epsilon=exp(\pi{}i/n)$ and $l\in(0,1,...,2n-1)$. 
The eigenvalues of $\textbf{M}$ are 

\begin{equation}
	 \lambda_j=\Sigma_{l=1}^{2n} M_{1 l}\epsilon^{j(l-1)}=m_0+2\Sigma_{k=1}^{n-1}m_k\cos(\pi{}jk/n)+(-1)^j m_n
 \end{equation} 
, and because of relation \ref{eq:NvsM}, those of $\textbf{S}$ are $\lambda_j^2$.
The eigenvector $\lambda_1$ is equal to the inflation factor 

 \begin{equation}
   L= m_0-m_n+2\Sigma_k^{\lfloor n/2 \rfloor} (m_k-m_{n-k}) cos(k\pi/n) 
 \end{equation}
\label{eq:InflationFactor}

A substitution tiling can only be a model set for a quasi crystal if its inflation factor is a Pisot- or PV-number \cite{Meyer95}, because a model set is point diffractive \cite{Hof95}. $L$ is a PV-number, if the absolute value of all its conjugates is less than 1. The conjugate eigenvalues are the ones for $j$ coprime to $n$. Using the above formulae we find that the inflation factors are PV-numbers in the following $m_1=1$ or $m_2=1$ or socalled \textit{single dent} cases:

\begin{table}[h]
	\centering
		\begin{tabular}{|l|l|l|l|l|}
		  \hline
			$m_0$ & $m_1$ & $m_2$ & $n$  \\ \hline 
			0  &  1 &   0 &   5    \\ 
			1  &  1 &   0 &   4,5,6,9   \\
			2  &  1 &   0 &   4,6 (Harriss)   \\
			1  &  0 &   1 &   5,7   \\
			2  &  0 &   1 &   5   \\
			\hline
		\end{tabular}
	\caption{Single dent substitution tiles with PV inflation factor.}
	\label{tab:PVInflationFactor}
\end{table}

The edge substitution matrix for a halfinteger edge sequence can be obtained by doubling the $n$ value. The fractional indices have to be doubled as well and become the $m_k$ values for odd k, whereas the $m_k$ for even k are zero. From table \ref{tab:PVInflationFactor} it is clear that the half integer single dent substitution tiles will not have PV inflation factors.

\section{Relationship with known rhomb tilings}
In all known rhomb tilings only positive area tiles are used. Nevertheless, some of the rules derived in this paper can also be applied. The examples we will discuss are the Ammann Beenker(AB) tiling, the Penrose (P) tiling, the binary or Lancon Billard (LB) tiling and the Harriss $n$-fold generalization of the Goodmann-Strauss tiling.  
The inflation factors of many known rhomb substitution tilings can be calculated using either formula \ref{eq:InflationFactor} or \ref{eq:NvsM}. One requirement is that the non-zero prototile edge angles occur in $\pm$ pairs, another that the area of the sum of the prototiles in the substition tile is equal to the substitution rhomb area. Than, the number of prototiles per substitution tile and, therefore, the substitution matrix can be calculated using equation \ref{eq:Nmatrix}, by neglecting the zero area tiles, and subtracting the negative tiles from the corresponding positive ones, or $n'_s=n_s-n_{2n-s}$.  This set of positive prototiles can be made even smaller by adding the ones having identical basic rhomb shapes or $n''_s=n_s+n_{n-s}-n_{n+s}-n_{2n-s}$.

The most famous examples are the AB ($n=4$) and P($n=5$) tilings in which the full set of 8, respectively 10 tiles have been reduced to a set of only 2 prototiles. The AB substitution tile has a (0, 1, -1) edge sequence. However, the edge configuration is different from the one in fig. \ref{fig:SubstitutionTileEdgeShape} and differs for both prototiles. The P substitution tiles have two different pairs of edges, a couple with a (1, -1)  and a couple with a (0, 2, -2) edge sequence. These edges have the same length, because $2\cos(\pi{}/5)$ is equal to $1+2\cos(2\pi{}/5)$. The corresponding sets of $m$-values, ($m_1=1$), and ($m_0=1, m_2=1$) respectively, also yield the same substitution matrix $\textbf{S}$ for $n=5$.

Some of the known rhomb tilings are in a way special cases of our general scheme. For instance, the binary or Lan\c{c}on Billard(LB) tiling. The LB substitution tiles have an $m_\frac{1}{2}=1$ edge. As was discussed in section \ref{Sec:FracMTiles}, the even $s$ tiles of the $n=5$ case are identical to the substitution tiles of the LB tiling, if substitution rule c) or d) is followed. Our treatment shows, that the LB tiling can be generalized to odd $n$ values different from five. The prototile set consists of all even $s$ tiles. The set of odd $s$ tiles produces another very similar tiling. Although the $s=n$ tile with a positive and negative part is involved, the prototiles in the substitution tiles are not cut up. For even $n$, the full set of prototiles, odd and even, are to be used to construct an LB-like tiling. 
 
Another special case of our tiling scheme is the Harriss tiling for arbitrary $n>3$ of which the Goodmann-Strauss tiling ($n=7$) is an example. The  edge sequence of the Harriss tiles is (0, 1, -1, 0), i.e. $m_0=2$ and $m_1=1$. The corresponding $n_t$ values are 6, 4 and 1 for $t=0$, 1 and 2 respectively. By replacing the $s\pm1$ tiles by $n-s\mp1$ tiles, as Harriss did, the set of prototiles is reduced to the even $s$ prototiles in the odd $n$ case. From the perspective of our scheme it has the advantage that the $s=n$ tile with a negative region is replaced by the $s=0$ tile whithout a negative region. For even $n$, Harriss had to introduce a number of additional tiles having extra dents or dimples. (Note, that the index $p$ used by Harris is equal to $n-s$)

Finally, the tiling with 11-fold rotation symmetry found by Maloney with computer assistence, is based on a set of only five complex prototiles with a (1,-1,-3,3,0,2,-2,-1,1,0,-5,5,-3,3,-1,1,4,-4,2,-2,0,-1,1,2,-2,-3,3,0,4,-4,-1,1,2,-2,0) edge sequence \cite{Maloney14}. The corresponding set of $(m_t)_0^5$ is (5,5,4,3,2,1). The positive and even $s$ tiles, are apparently sufficient to tile the plane.

\section{Gaps and Overlaps}
A consequence of the admittance of negative prototiles in our general model is that members of the substitution tile set might contain negative parts, gaps and overlaps. In the following we will analyze the conditions under which these defects disappear in tilings.

\begin{figure}[!ht]
	\centering
		\includegraphics[width=0.8\textwidth]{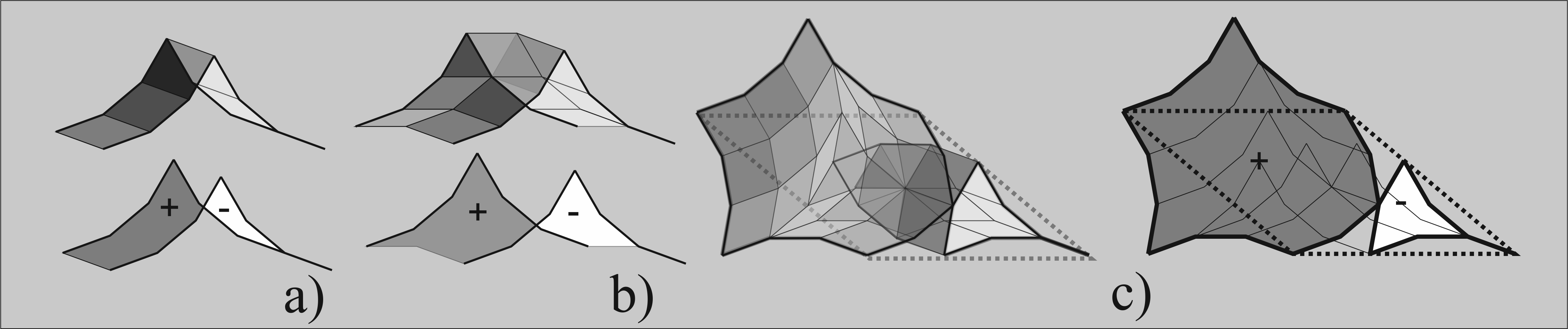}
	\caption{\small {Tile construction in terms of rows of prototiles or worms.
a) Row of prototiles along an edge. At the crossing point of the upper and lower worm edges the exterior negative and positive prototile parts compensate each other and the sign of the interior prototile (parts) changes. b) By adding the next row of prototiles one gets a double worm following the same rule, i.e. the interior parts are either negative or positive and the exterior is empty. c) Once the tile is completed, it will consists of positive and negative regions, if its edges cross. The example shown here is a $T_7^9$ tile with a (1,2,3,-3,-2,-1) edge sequence.}}
	\label{fig:worms}
\end{figure}

As was suggested earlier, the general substitution rule for an $N\times N$ substitution tile may be formulated in terms of worms of $N$ prototiles. The procedure is illustrated in figure \ref{fig:worms}.  The worms are bounded by two parallel copies of a substitution tile edge and separated by one edge unit along neighboring edges. The interior of the worm, i.e. the space between the two edges, is either filled with positive or negative tiles or tile pieces, provided the edge does not contain loops. The sign changes at the crossing points of the two. Adding a neighboring worm, one gets a double worm with congruent edges, and interior regions which are either positive or negative.  If the tile edges do not cross anywhere after completion of the tile, i.e. if the circumpherence does not contain loops,  the filling of the tile is a single positive layer of tiles or tile pieces. If they do cross somewhere, there will be regions filled with a layer of negative tiles or tile pieces. To fully tile the plane we may repeatedly inflate a positive tile and apply the edge substitution rule. This involves the substitution of the inflated edge segments by the initial edge sequence. Since the edge substitution may be done in two orientations related by a twofold rotation, each time $2^N$ different next generation edge sequences are possible. The new extended edge contains $N^g$ segments, $g$ being the substitution generation. Consequently, the worms consisting of $N^g$ prototiles are the new building blocks of the tile. If we do not admit loops in the edge sequence, the interior of the multiply inflated tile will consist of positive and negative regions separated by the crossing points of the four edges, but no gaps or overlaps.

\section{Summary}
A general scheme for constructing a non-periodic rhomb tiling was derived for arbitrary $n\geq 3$. It admits zero and negative area tiles to ensure generality. If the edge sequence does not contain loops, such a tiling will consist of single layered positive and possibly negative regions, but there will be no gaps or overlaps in the tiling. Overlaps occurring in the tiling process are always annihilated by negative or subtraction tiles. This generally, but not always, leads to the cutting up of tiles.

The corresponding edge substitution matrix can be expressed as a circulant matrix $\textbf{M}$, the inflation factor being one of its eigenvalues. The tile substitution matrix $\textbf{S}$ is equal to the square of $\textbf{M}$. For a given edge sequence this equality determines  the type and number of prototiles in a substitution tile. A possible arrangement of the prototiles is obtained by drawing a grid using parallel copies of the edge. This arrangement can be expressed in the form of a matrix. 
   
The smallest ($2\times2$) substitution tiles were discussed in detail. Two different types can be distingushed, those with an integer and those with a halfinteger edge sequence. The latter are of special interest, because some of the tilings do not contain tiles cut up by the application of subtraction tiles. An example is the LB tiling ($n=5$), which can be generalized to arbitrary $n$.    

\section*{Acknowledgement}
I gratefully acknowledge the very useful discussions with Dr. D. Frettl\"{o}h from the university of Bielefeld.

\end{document}